\documentclass[preprint,11pt]{article}

\usepackage{graphicx,epsf,amsbsy,amsmath,amssymb}
\usepackage[latin1]{inputenc}

\newtheorem{prop}{Proposition}[section] % Les propositions sont numérotées comme

\begin{document}

\title{Slow Manifold of a Neuronal Bursting Model}

\author{Jean-Marc Ginoux and Bruno Rossetto, \\ P.R.O.T.E.E. Laboratory, I.U.T. de Toulon,\\ Universit\'e du Sud,
B.P. 20132, 83957, La Garde cedex, France,\\ E-mail: ginoux@univ-tln.fr, rossetto@univ-tln.fr}

\date{{\bf Keywords}: differential geometry; curvature;
torsion; slow-fast dynamics; neuronal bursting models.}

\maketitle

\begin{abstract}
Comparing neuronal bursting models (NBM) with slow-fast autonomous
dynamical systems (S-FADS), it appears that the specific features of
a (NBM) do not allow a determination of the analytical slow
manifold equation with the singular approximation method. So, a new approach based on Differential
Geometry, generally used for (S-FADS), is proposed. Adapted to (NBM), this new method provides three
equivalent manners of determination of the analytical slow manifold
equation. Application is made for the three-variables model of
neuronal bursting elaborated by Hindmarsh and Rose which is one of
the most used mathematical representation of the widespread
phenomenon of oscillatory burst discharges that occur in real
neuronal cells.
\end{abstract}

\section{Slow-fast autonomous dynamical systems, neuronal bursting models}
\label{sec1}

\subsection{Dynamical systems}

In the following we consider a system of differential equations
defined in a compact E included in $\mathbb{R}$:

\begin{equation}
\label{eq1}
\frac{d{\vec {X}}}{dt} = \vec \Im \left( \vec {X} \right)
\end{equation}

with

\[
\vec {X} = \left[ {x_1 ,x_2 ,...,x_n } \right]^t \in E \subset
\mathbb{R}^n
\]

and

\[
\vec \Im \left( \vec {X} \right) = \left[ {f_1 \left( \vec {X}
\right),f_2 \left( \vec {X} \right),...,f_n \left( \vec {X} \right)}
\right]^t \in E \subset \mathbb{R}^n
\]

The vector $\vec \Im \left( \vec {X} \right)$ defines a velocity
vector field in E whose components $f_i $ which are supposed to be
continuous and infinitely differentiable with respect to all $x_i $
and $t$, i.e., are $C^\infty $ functions in E and with values
included in $\mathbb{R}$, satisfy the assumptions of the
Cauchy-Lipschitz theorem. For more details, see for example
\cite{2}. A solution of this system is an integral curve $\vec
{X}\left( t \right)$ tangent to $\vec \Im $ whose values define the
\textit{states} of the \textit{dynamical system} described by the Eq. (\ref{eq1}). Since none of the
components $f_i $ of the velocity vector field depends here
explicitly on time, the system is said to be \textit{autonomous}.

\subsection{Slow-fast autonomous dynamical system (S-FADS)}

A (S-FADS) is a \textit{dynamical system} defined under the same conditions as
previously but comprising a small multiplicative parameter
$\varepsilon$ in one or several components of its velocity vector
field:

\begin{equation}
\label{eq2} \frac{d{\vec {X}}}{dt} = \vec \Im_\varepsilon
\left( \vec {X} \right)
\end{equation}

with

\[
\vec \Im_\varepsilon \left( \vec {X} \right) = \left[
{\frac{1}{\varepsilon}f_1 \left( \vec {X} \right),f_2 \left( \vec
{X} \right),...,f_n \left( \vec {X} \right)} \right]^t \in E \subset
\mathbb{R}^n
\]

\[
0 < \varepsilon \ll 1
\]

The functional jacobian of a (S-FADS) defined by (\ref{eq2}) has an eigenvalue called
``fast'', i.e., great on a large domain of the phase space. Thus, a
``fast'' eigenvalue is expressed like a polynomial of valuation $ -
1$ in $\varepsilon $ and the eigenmode
which is associated with this ``fast'' eigenvalue is said:\\

- ``evanescent'' if it is negative,

- ``dominant'' if it is positive.\\

The other eigenvalues called ``slow'' are expressed like a
polynomial of valuation $0$ in $\varepsilon $.

\newpage

\subsection{Neuronal bursting models (NBM)}

A (NBM) is a \textit{dynamical system} defined under the same
conditions as previously but comprising a large multiplicative
parameter $\varepsilon^{ - 1}$ in one component of its velocity
vector field:

\begin{equation}
\label{eq3} \frac{d{\vec {X}}}{dt} = \vec \Im_\varepsilon
\left( \vec {X} \right)
\end{equation}

with

\[
\vec \Im_\varepsilon \left( \vec {X} \right) = \left[ { f_1 \left(
\vec {X} \right),f_2 \left( \vec {X} \right),...,\varepsilon f_n
\left( \vec {X} \right)} \right]^t \in E \subset \mathbb{R}^n
\]

\[
0 < \varepsilon \ll 1
\]

The presence of the multiplicative parameter $\varepsilon^{ - 1}$ in
one of the components of the velocity vector field makes it possible
to consider the system (\ref{eq3}) as a kind of
\textit{slow}-\textit{fast} autonomous dynamical system
(S-FADS). So, it possesses a \textit{slow manifold}, the
equation of which may be determined. But, paradoxically, this model
is not \textit{slow-fast} in the sense defined previously. A
comparison between three-dimensional (S-FADS) and
(NBM) presented in Table 1 emphasizes their differences. The dot $(\cdot)$
represents the derivative with respect to time and $\varepsilon \ll
1$.

\begin{table}[htbp]
\centering
\caption[Comparison between slow-fast autonomous
dynamical systems and neuronal bursting models]{Comparison between
(S-FADS) and (NBM)}

\vspace{0.1in}

\begin{tabular}
{|c|c|c|}
\hline \multicolumn{2}{|c|}{\textbf{(S-FADS) vs (NBM)} \par  }  \\
\hline
\mbox{   } \par & \mbox{   } \par \\
$ \hspace{0.1in} \dfrac{d{\vec {X}}}{dt} \left(
{{\begin{array}{*{20}c}
 \dot {x} \hfill \\
 \dot {y} \hfill \\
 \dot {z} \hfill \\
\end{array} }} \right) = \vec{\Im_\varepsilon} \left( {{\begin{array}{*{20}c}
 {\dfrac{1}{\varepsilon } f\left( {x,y,z} \right)} \vspace{4pt} \\
 {\mbox{ }g\left( {x,y,z} \right)} \vspace{4pt} \\
 {h\left( {x,y,z} \right)} \\
\end{array} }} \right) \mbox{ }$ \par & $ \hspace{0.1in} \dfrac{d{\vec {X}}}{dt} \left( {{\begin{array}{*{20}c}
 \dot {x} \hfill \\
 \dot {y} \hfill \\
 \dot {z} \hfill \\
\end{array} }} \right) = \vec{\Im_\varepsilon} \left( {{\begin{array}{*{20}c}
 {\mbox{ }f\left( {x,y,z} \right)} \vspace{4pt} \\
 {\mbox{ }g\left( {x,y,z} \right)} \vspace{4pt} \\
 {\varepsilon h\left( {x,y,z} \right)} \\
\end{array} }} \right) \mbox{ } $ \par  \\
\mbox{   } \par & \mbox{   } \par \\
\hline
\mbox{   } \par & \mbox{   } \par \\
$ \hspace{0.1in} \dfrac{d{\vec {X}}}{dt} \left(
{{\begin{array}{*{20}c}
\dot {x} \hfill \\
\dot {y} \hfill \\
\dot {z} \hfill \\
\end{array} }} \right) = \vec{\Im_\varepsilon} \left( {{\begin{array}{*{20}c}
 {\mbox{ } fast} \vspace{4pt} \\
 {\mbox{ }slow} \vspace{4pt} \\
 {slow} \\
\end{array} }} \right) \mbox{ } $ \par & $\hspace{0.1in} \dfrac{d{\vec {X}}}{dt} \left( {{\begin{array}{*{20}c}
 \dot {x} \hfill \\
 \dot {y} \hfill \\
 \dot {z} \hfill \\
\end{array} }} \right) = \vec{\Im_\varepsilon} \left( {{\begin{array}{*{20}c}
 {\mbox{ }fast} \vspace{4pt} \\
 {\mbox{ }fast} \vspace{4pt} \\
 {slow} \\
\end{array} }} \right) \mbox{ }$ \par  \\
\mbox{   } \par & \mbox{   } \par \\
 \hline
\end{tabular}
\label{tab1}
\end{table}

\newpage

\section{Analytical slow manifold equation}
\label{sec2}

There are many methods of determination of the analytical equation
of the slow manifold. The classical one based on the singular
perturbations theory \cite{1} is the so-called \textit{singular
approximation method}. But, in the specific case of a
(NBM), one of the hypothesis of the
Tihonov's theorem is not checked since the \textit{fast dynamics} of
the \textit{singular approximation}
has a periodic solution. Thus, another approach developed in
\cite{4} which consist in using \textit{Differential
Geometry} formalism may be used.

\subsection{Singular approximation method}

The\textit{ singular approximation} of the \textit{fast} dynamics
constitutes a quite good approach since the third component of the
velocity is very weak and so, $z$ is nearly constant along the
periodic solution. In dimension three the system (\ref{eq3}) can be
written as a system of differential equations defined in a compact E
included in $\mathbb{R}$:

\[
\hspace{0.1in} \frac{d{\vec {X}}}{dt} = \left(
{{\begin{array}{*{20}c}
 {\frac{d{x}}{dt}} \vspace{4pt} \\
 {\frac{d{y}}{dt}} \vspace{4pt} \\
 {\frac{d{z}}{dt}} \\
\end{array} }} \right) = \vec \Im_\varepsilon \left( {{\begin{array}{*{20}c}
 {\mbox{ }f\left( {x,y,z} \right)} \vspace{4pt} \\
 {\mbox{ }g\left( {x,y,z} \right)} \vspace{4pt} \\
 {\varepsilon h\left( {x,y,z} \right)} \\
\end{array} }} \right)
\]

On the one hand, since the system (\ref{eq3}) can be considered as a
(S-FADS), the
\textit{slow} dynamics of the \textit{singular
approximation} is given by:

\begin{equation}
\label{eq4} \hspace{0.4in} \left\{ {{\begin{array}{*{20}c}
 {f\left( {x,y,z} \right) = 0} \hfill \\
 {g\left( {x,y,z} \right) = 0} \hfill \\
\end{array} }} \right.
\end{equation}

The resolution of this reduced system composed of the two first
equations of the right hand side of (\ref{eq3}) provides a
one-dimensional \textit{singular manifold}, called \textit{singular
curve}. This curve doesn't play any role in the construction of the
periodic solution. But we will see that there exists all the more a
\textit{slow dynamics}. On the other hands, it presents a
\textit{fast} dynamics which can be given while posing the following
change:

\[
\tau = \varepsilon t\mbox{ } \Leftrightarrow \mbox{
}\frac{d}{dt} = \varepsilon \frac{d}{d{\tau} }
\]

\newpage

The system (\ref{eq3}) may be re-written as:

\begin{equation}
\label{eq5}
\frac{d{\vec {X}}}{d{\tau} } = \left(
{{\begin{array}{*{20}c}
 {\frac{d{x}}{d{\tau} }} \vspace{4pt} \\
 {\frac{d{y}}{d{\tau} }} \vspace{4pt} \\
 {\frac{d{z}}{d{\tau} }} \\
\end{array} }} \right) = \vec{\Im_\varepsilon} \left( {{\begin{array}{*{20}c}
 {\varepsilon ^{ - 1}f\left( {x,y,z} \right)} \vspace{4pt} \\
 {\varepsilon ^{ - 1}g\left( {x,y,z} \right)} \vspace{4pt} \\
 {\mbox{ }h\left( {x,y,z} \right)} \\
\end{array} }} \right)
\end{equation}

So, the \textit{fast} dynamics of the \textit{singular
approximation} is provided by the study of the reduced system
composed of the two first equations of the right hand side of
(\ref{eq5}).

\begin{equation}
\label{eq6}
\left. {\frac{d{\vec {X}}}{d{\tau} }} \right|_{fast} =
\left( {{\begin{array}{*{20}c}
 {\frac{d{x}}{d{\tau} }} \vspace{4pt} \\
 {\frac{d{y}}{d{\tau} }}  \\
\end{array} }} \right) = \vec{\Im_\varepsilon} \left( {{\begin{array}{*{20}c}
 {\varepsilon ^{ - 1}f\left( {x,y,z^\ast } \right)} \vspace{4pt} \\
 {\varepsilon ^{ - 1}g\left( {x,y,z^\ast } \right)} \\
\end{array} }} \right)
\end{equation}

Each point of the \textit{singular curve} is a singular point of the
\textit{singular approximation} of the
\textit{fast} dynamics. For the $z$ value for which there is a
periodic solution, the \textit{singular approximation} exhibits an
unstable focus, attractive with respect to the \textit{slow}
eigendirection.

\subsection{Differential Geometry formalism}

Now let us consider a three-dimensional system defined by
(\ref{eq3}) and let's define the instantaneous acceleration vector
of the \textit{trajectory curve }$\vec {X}\left( t \right)$. Since
the functions $f_i $ are supposed to be $C^\infty$ functions in a
compact E included in $\mathbb{R}$, it is possible to calculate the
total derivative of the vector field $\vec{\Im_\varepsilon}$. As the
instantaneous vector function $\vec V \left( t \right)$ of the
scalar variable $t$ represents the velocity vector of the mobile M
at the instant $t$, the total derivative of $\vec V \left( t
\right)$ is the vector function $\vec \gamma \left( t \right)$ of
the scalar variable $t$ which represents the instantaneous
acceleration vector of the mobile M at the instant $t$. It is noted:

\begin{equation}
\label{eq7} \vec \gamma \left( t \right)\mbox{ } = \mbox{
}\frac{d{\vec V} \left( t \right)}{dt}
\end{equation}

Even if neuronal bursting models are not exactly slow-fast
autonomous dynamical systems, the new approach of determining the
\textit{slow manifold} equation developed in
\cite{4} may still be applied. This method is using
\textit{Differential Geometry}
properties such as \textit{curvature} and
\textit{torsion} of the \textit{trajectory curve}
$\vec {X}\left( t \right)$, integral of \textit{dynamical systems}
to provide their \textit{slow manifold} equation.

\newpage

\begin{prop}
\label{prop un}

The location of the points where the local torsion of the trajectory
curves integral of a dynamical system defined by (\ref{eq3})
vanishes, provides the analytical equation of the slow manifold
associated with this system.

\begin{equation}
\label{eq8} \frac{1}{\Im } = - \frac{\dot {\vec {\gamma }} \cdot
\left( {\vec \gamma \times \vec V } \right)}{\left\| {\vec \gamma
\times \vec V } \right\|^2} = 0\mbox{ } \Leftrightarrow \mbox{ }\dot
{\vec {\gamma }} \cdot \left( {\vec \gamma \times \vec V } \right) =
0
\end{equation}

Thus, this equation represents the \textit{slow manifold} of a neuronal bursting model defined by (\ref{eq3}).

\end{prop}

The particular features of neuronal bursting models (\ref{eq3}) will
lead to a simplification of this Proposition 1. Due to
the presence of the small multiplicative parameter $\varepsilon $ in
the third components of its velocity vector field, instantaneous
velocity vector $\vec V \left( t \right)$ and instantaneous
acceleration vector $\vec \gamma \left( t \right)$ of the model
(\ref{eq3}) may be written:

\begin{equation}
\label{eq9} \hspace{0.4in} \vec V \left( {{\begin{array}{*{20}c}
 \dot {x} \hfill \\
 \dot {y} \hfill \\
 \dot {z} \hfill \\
\end{array} }} \right) = \vec \Im_\varepsilon \left( {{\begin{array}{*{20}c}
 {O\left( {\varepsilon ^0} \right)} \vspace{4pt} \\
 {O\left( {\varepsilon ^0} \right)} \vspace{4pt} \\
 {O\left( {\varepsilon ^1} \right)} \vspace{4pt} \\
\end{array} }} \right)
\end{equation}

and

\begin{equation}
\label{eq10}
\hspace{0.3in} \vec {\gamma }\left(
{{\begin{array}{*{20}c}
 \ddot {x} \hfill \\
 \ddot {y} \hfill \\
 \ddot {z} \hfill \\
\end{array} }} \right) = \frac{d{\vec \Im_\varepsilon }}{dt}\left(
{{\begin{array}{*{20}c}
 {O\left( {\varepsilon ^1} \right)} \vspace{4pt} \\
 {O\left( {\varepsilon ^1} \right)} \vspace{4pt} \\
 {O\left( {\varepsilon ^2} \right)} \vspace{4pt} \\
\end{array} }} \right)
\end{equation}

where ${O\left( {\varepsilon ^n} \right)}$ is a polynomial of $n$
degree in $\varepsilon$

Then, it is possible to express the vector product $\vec V \times
\vec \gamma $ as:

\begin{equation}
\label{eq11} \hspace{0.4in} \vec V \times \vec \gamma = \left(
{{\begin{array}{*{20}c}
 {\dot {y}\ddot {z} - \ddot {y}\dot {z}} \hfill \\
 {\ddot {x}\dot {z} - \dot {x}\ddot {z}} \hfill \\
 {\dot {x}\ddot {y} - \ddot {x}\dot {y}} \hfill \\
\end{array} }} \right)
\end{equation}

Taking into account what precedes (\ref{eq9}, \ref{eq10}), it
follows that:

\begin{equation}
\label{eq12} \hspace{0.4in}\vec V \times \vec \gamma = \left(
{{\begin{array}{*{20}c}
 {O\left( {\varepsilon ^2} \right)} \vspace{4pt} \\
 {O\left( {\varepsilon ^2} \right)} \vspace{4pt} \\
 {O\left( {\varepsilon ^1} \right)} \vspace{4pt} \\
\end{array} }} \right)
\end{equation}

So, it is obvious that since $\varepsilon $ is a small parameter,
this vector product may be written:

\begin{equation}
\label{eq13} \hspace{0.4in}\vec V \times \vec \gamma \approx \left(
{{\begin{array}{*{20}c}
 0 \\
 0 \\
 {O\left( {\varepsilon ^1} \right)} \\
\end{array} }} \right)
\end{equation}

Then, it appears that if the third component of this vector product
vanishes when both instantaneous velocity vector $\vec V \left( t
\right)$ and instantaneous acceleration vector $\vec \gamma \left( t
\right)$ are collinear. This result is particular to this kind of
model which presents a small multiplicative parameter in one of the
right-hand-side component of the velocity vector field and makes it
possible to simplify the previous Proposition 1.\\

\begin{prop}
\label{prop deux}

The location of the points where the instantaneous velocity vector
$\vec V \left( t \right)$ and instantaneous acceleration vector
$\vec \gamma \left( t \right)$ of a neuronal bursting model defined
by (\ref{eq3}) are collinear provides the analytical equation of the
\textit{slow manifold} associated with this
dynamical system.

\begin{equation}
\label{eq14} \vec V \times \vec \gamma = \vec 0  \mbox{ }
\Leftrightarrow \mbox{ }\dot {x}\ddot {y} - \ddot {x}\dot {y} = 0
\end{equation}

\end{prop}

Another method of determining the \textit{slow manifold} equation proposed in \cite{14} consists in considering the
so-called \textit{tangent linear system approximation}. Then, a coplanarity condition between
the instantaneous velocity vector $\vec V \left( t \right)$ and the
\textit{slow }eigenvectors of the \textit{tangent linear system}
gives the \textit{slow manifold} equation.

\begin{equation}
\label{eq15} \vec V .\left( {\vec {Y_{\lambda _2 } } \times \vec
{Y_{\lambda _3 } } } \right) = 0
\end{equation}

where $\vec {Y_{\lambda _i } } $ represent the \textit{slow}
eigenvectors of the \textit{tangent linear system}. But, if these
eigenvectors are complex the \textit{slow manifold} plot may be interrupted. So, in order to avoid such
inconvenience, this equation has been multiplied by two conjugate
equations obtained by circular permutations.

\[
\left[ {\vec V \cdot \left( {\vec {Y_{\lambda _2 } } \times \vec
{Y_{\lambda _3 } } } \right)} \right] \cdot \left[ {\vec V \cdot
\left( {\vec {Y_{\lambda _1 } } \times \vec {Y_{\lambda _2 } } }
\right)} \right] \cdot \left[ {\vec V \cdot \left( {\vec {Y_{\lambda
_1 } } \times \vec {Y_{\lambda _3 } } } \right)} \right] = 0
\]

It has been established in \cite{4} that this real analytical
\textit{slow manifold} equation can be written:

\begin{equation}
\label{eq16} \left( {J^2\vec V } \right) \cdot \left( {\vec \gamma
\times \vec V } \right) = 0
\end{equation}

since the the \textit{tangent linear system
approximation} method
implies to suppose that the functional jacobian matrix is
stationary. That is to say

\[
\frac{d J}{dt} = 0
\]

and so,

\[ \dot {\vec {\gamma }} = J\frac{d{\vec V }}{d t} +
\frac{d J}{d t}\vec V = J\vec \gamma + \frac{d J}{d t}\vec V =
J^2\vec V + \frac{d J}{d t}\vec V \approx J^2\vec V
\]

\vspace{0.1in}

\begin{prop}
\label{prop trois}

The coplanarity condition (\ref{eq15}) between the instantaneous
velocity vector and the slow eigenvectors of the tangent linear
system transformed into the real analytical equation (\ref{eq16})
provides the \textit{slow manifold} equation of
a neuronal bursting model defined by (\ref{eq3}).
\end{prop}

\section{Application to a neuronal bursting model}
\label{sec3}

The transmission of nervous impulse is secured in the brain by
action potentials. Their generation and their rhythmic behaviour are
linked to the opening and closing of selected classes of ionic
channels. The membrane potential of neurons can be modified by
acting on a combination of different ionic mechanisms. Starting from
the seminal works of Hodgkin-Huxley \cite{7,11} and FitzHugh-Nagumo
\cite{3,12}, the Hindmarsh-Rose \cite{6,13} model consists of three
variables: $x$, the membrane potential, $y$, an intrinsic current
and $z$, a \textit{slow} adaptation current.

\subsection{Hindmarsh-Rose model of bursting neurons}
% Always give a unique label
% and use \ref{<label>} for cross-references
% and \cite{<label>} for bibliographic references
% use \sectionmark{}
% to alter or adjust the section heading in the running head

\begin{equation}
\label{eq17}
\hspace{0.4in}\left\{ {{\begin{array}{{ll}}
 {\frac{d{x}}{dt} = y - f\left( x \right) - z + I} \vspace{4pt} \hfill \\
 {\frac{d{y}}{dt} = g\left( x \right) - y}  \vspace{4pt} \hfill \\
 {\frac{d{z}}{dt} = \varepsilon \left( {h\left( x \right) - z} \right)} \hfill \\
\end{array}}} \right.
\end{equation}

$I$ represents the applied current, $f\left( x \right) = ax^3 -
bx^2$ and $g\left( x \right) = c - dx^2$ are respectively cubic and
quadratic functions which have been experimentally deduced \cite{5}.
$\varepsilon$ is the time scale of the slow adaptation current and
$h\left( x \right) = x - x^\ast$ is the scale of the influence of
the \textit{slow} dynamics, which determines whether the neuron
fires in a tonic or in a burst mode when it is exposed to a
sustained current input and where $\left( {x^\ast ,y^\ast } \right)$
are the co-ordinates of the leftmost equilibrium point of the model
(1) without adaptation, i.e., $I = 0$.\\ Parameters used for
numerical simulations are:\\ $a = 1$, $b = 3$, $c = 1$, $d = 5$,
$\varepsilon = 0.005$, $s = 4$, $x^\ast = \frac{ - 1 - \sqrt 5 }{2}$
and $I =
3.25$.\\

While using the method proposed in the section 2 it is
possible to determine the analytical \textit{slow
manifold} equation of the Hindmarsh-Rose
84'model \cite{6}.

\subsection{Slow manifold of the Hindmarsh-Rose 84'model}

In Fig. 1 is presented the \textit{slow manifold} of the Hindmarsh-Rose 84'model determined with the Proposition 1.\\

\begin{figure}[htbp]
\centering
\includegraphics[width=10cm,height=10cm]{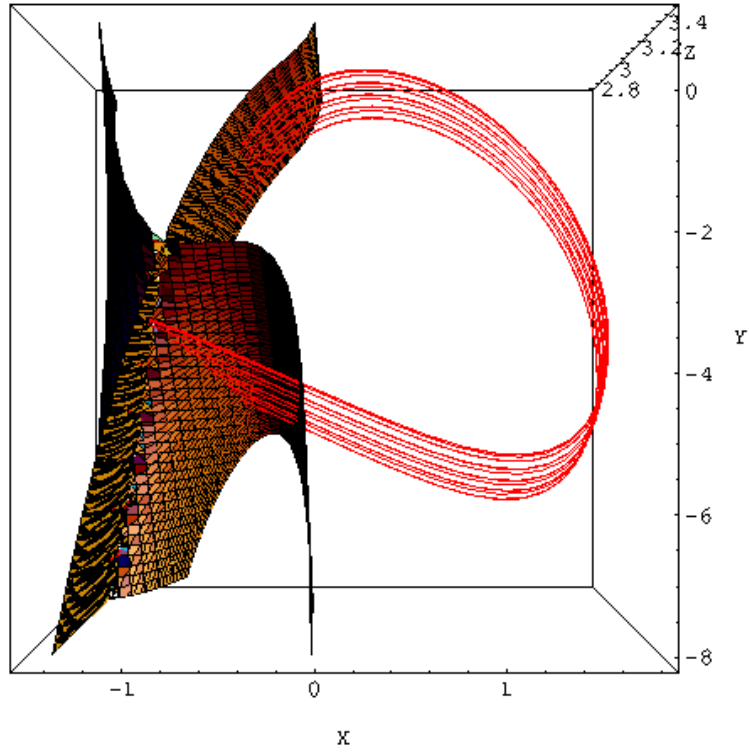}
\caption{\textit{Slow manifold} of the Hindmarsh-Rose 84'model with the Prop. 1.}
\label{fig:1}
\end{figure}

\newpage

The \textit{slow manifold} provided with the use of the \textit{collinearity condition} between both instantaneous velocity vector and instantaneous
acceleration vector, i.e., while using the
Proposition 2 is presented in Fig. 2.\\

\begin{figure}[htbp]
\centering
\includegraphics[width=10cm,height=10cm]{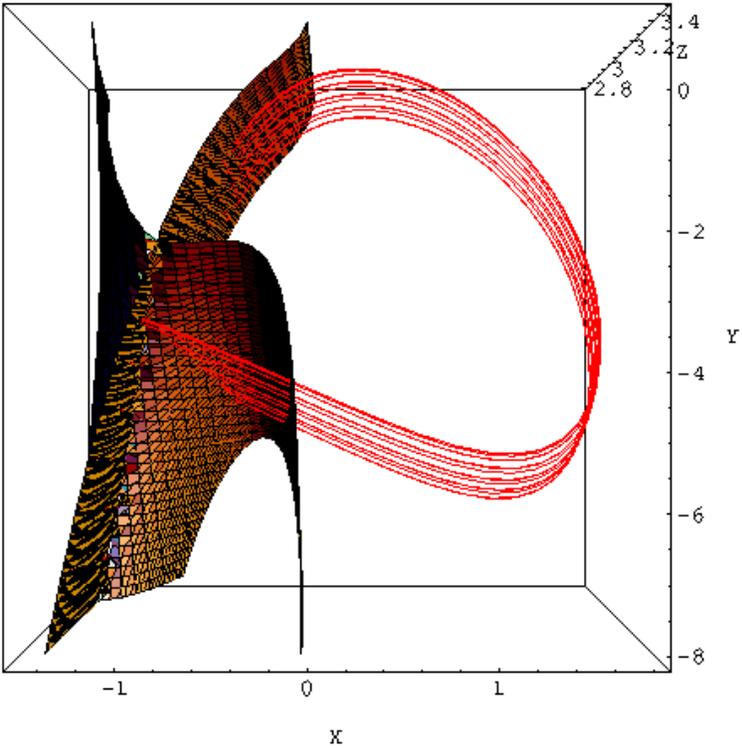}
\caption{\textit{Slow manifold} of the Hindmarsh-Rose 84'model with the Prop. 2.}
\label{fig:2}
\end{figure}

\newpage

Figure 3 presents the \textit{slow manifold} of the Hindmarsh-Rose
84'model obtained with the \textit{tangent linear system
approximation}, i.e.,
with the use of Proposition 3.

\begin{figure}[htbp]
\centering
\includegraphics[width=10cm,height=10cm]{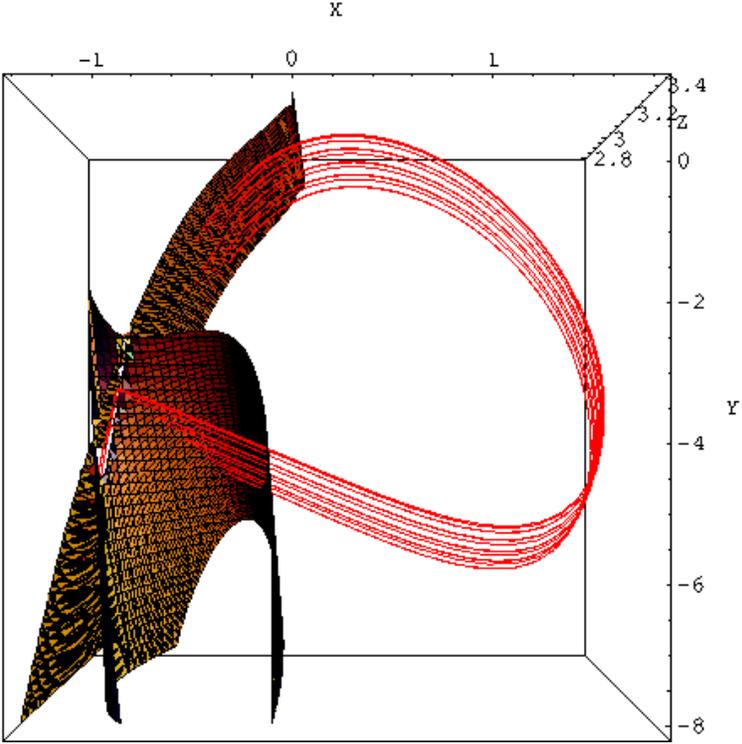}
\caption{\textit{Slow manifold} of the Hindmarsh-Rose 84'model with the Prop. 3.}
\label{fig:3}
\end{figure}

\newpage

\section{Discussion}
\label{sec4}

Since in the case of neuronal bursting model (NBM) bursting models one of the Tihonov's hypothesis is not checked, the classical \textit{singular approximation method} can not be used to determine the analytical \textit{slow manifold} equation. In this work the application of the \textit{Differential Geometry} formalism provides new alternative methods of determination of the \textit{slow manifold} equation of a neuronal bursting model (NBM).

\begin{itemize}

\item the \textit{torsion method}, i.e., the location of the points where the local
\textit{torsion} of the \textit{trajectory curve}, integral of
\textit{dynamical systems} vanishes,

\item the \textit{collinearity condition} between the instantaneous
velocity vector $\overrightarrow V $, the instantaneous acceleration
vector $\overrightarrow \gamma $,

\item the \textit{tangent linear system approximation}, i.e., the coplanarity condition
between the instantaneous velocity vector eigenvectors transformed
into a real analytical equation.

\end{itemize}

The striking similarity of all figures due to the smallness of the
parameter $\varepsilon$ highlights the equivalence between all the
propositions. Moreover, even if the presence of this small parameter
$\varepsilon$ in one of the right-hand-side component of the
instantaneous velocity vector field of a (NBM) prevents from using the \textit{singular
approximation method}, it clarifies the Proposition 1
and transforms it into a \textit{collinearity condition} in
dimension three, i.e., Proposition 2. Comparing
(S-FADS) and (NBM) in Table 1 it can be noted that in a (S-FADS) there is one
fast component and two fast while in a (NBM) the situation is
exactly reversed. Two fast components and one slow. So, considering
(NBM) as a particular class of (S-FADS) we suggest to call (NBM)
\textit{fast-slow} instead of \textit{slow-fast} in order to avoid
any confusion. Further research should highlight other specific
features of (NBM).

\section*{Acknowledgements}

Authors would like to thank Professors M. Aziz-Alaoui and
C. Bertelle for their useful collaboration.


\begin{thebibliography}{08}

\bibitem{1} Andronov AA, Khaikin SE, {\&} Vitt AA (1966)
Theory of oscillators, Pergamon Press, Oxford


\bibitem{2} Coddington EA {\&} Levinson N, (1955)
Theory of Ordinary Differential Equations, Mac Graw Hill, New York


\bibitem{3} Fitzhugh R (1961)
Biophys. J 1:445--466


\bibitem{4} Ginoux JM {\&} Rossetto B (2006)
Int. J. Bifurcations and Chaos, (in print)


\bibitem{5} Hindmarsh JL {\&} Rose RM (1982) Nature 296:162--164


\bibitem{6} Hindmarsh JL {\&} Rose RM (1984)
Philos. Trans. Roy. Soc. London Ser. B 221:87--102


\bibitem{7} Hodgkin AL {\&} Huxley AF (1952)
J. Physiol. (Lond.) 116:473--96



\bibitem{8} Hodgkin AL {\&} Huxley AF (1952)
J. Physiol. (Lond.) 116: 449--72



\bibitem{9} Hodgkin AL {\&} Huxley AF (1952)
J. Physiol. (Lond.) 116: 497--506



\bibitem{10} Hodgkin AL {\&} Huxley AF (1952)
J. Physiol. (Lond.) 117: 500--44


\bibitem{11} Hodgkin AL, Huxley AF {\&} Katz B (1952) B. Katz
J. Physiol. (Lond.) 116: 424--48


\bibitem{12} Nagumo JS, Arimoto S {\&} Yoshizawa S (1962)
Proc. Inst. Radio Engineers 50:2061--2070


\bibitem{13} Rose RM {\&} Hindmarsh JL (1985)
Proc. R. Soc. Ser. B 225:161--193

\bibitem{14} Rossetto B, Lenzini T, Suchey G  {\&} Ramdani S (1998)
Int. J. Bifurcation and Chaos, vol. 8 (11):2135-2145

\end{thebibliography}
\end{document}